\documentclass[11pt]{article}
\usepackage{amsfonts}

\newtheorem{Theorem}{\sc Theorem}[section]
\newtheorem{Proposition}[Theorem]{\sc Proposition}
\newtheorem{Corollary}[Theorem]{\sc Corollary}
\newtheorem{Definition}[Theorem]{\sc Definition}
\newtheorem{Construction}[Theorem]{\sc Construction}
\newtheorem{Lemma}[Theorem]{\sc Lemma}
\newtheorem{Remark}[Theorem]{\sc Remark}
\input{tcilatex}

\begin{document}

\title{A RADON-NIKODYM THEOREM\ FOR\ COMPLETELY\ MULTI-POSITIVE LINEAR\ MAPS\ AND\
ITS\ APPLICATIONS }
\author{MARIA\ JOI\c{T}A}
\maketitle

\begin{abstract}
A completely $n$ -positive linear map from a locally $C^{*}$-algebra $A$ to
another locally $C^{*}$-algebra $B\;$is an $n\times n$ matrix whose elements
are continuous linear maps from $A$ to $B$ and which verifies the condition
of completely positivity. In this paper we prove a Radon-Nikodym type
theorem for strict completely $n$-positive linear maps which describes the
order relation on the set of all strict completely $n$ -positive linear maps
from a locally $C^{*}$-algebra $A$ to a $C^{*}$-algebra $B$, in terms of a
self-dual Hilbert $C^{*}$-module structure induced by each strict completely 
$n$ -positive linear map. As applications of this result we characterize the
pure completely $n$-positive linear maps from $A$ to $B$ and the extreme
elements in the set of all identity preserving completely $n$-positive
linear maps from $A$ to $B$. Also we determine a certain class of extreme
elements in the set of all identity preserving completely positive linear
maps from $A$ to $M_n(B)$.

\end{abstract}

\section{Introduction and preliminaries}

The concept of matricial order plays an important role to understand the
infinite -dimensional non-commutative structure of operator algebras.
Completely positive maps, as the natural ordering attached to this structure
have been studied extensively [\textbf{1}, \textbf{2}, \textbf{5}, \textbf{4}%
, \textbf{7}, \textbf{8, 9}, \textbf{10}, \textbf{11}, \textbf{15}, \textbf{%
16}, \textbf{18}, \textbf{19}].

Given a $C^{*}$-algebra $A$ we denote by $M_n(A)$ the $C^{*}$-algebra of all 
$n\times n$ matrices with elements in $A.$

\begin{Definition}
A completely positive map from a $C^{*}$-algebra $A$ to another $C^{*}$%
-algebra $B$ is a linear map $\rho :A\rightarrow B$ such that the linear map 
$\rho _n:M_n(A)\rightarrow M_n(B)$ defined by 
\[
\rho _n\left( \left[ a_{ij}\right] _{i,j=1}^n\right) =\left[ \rho \left(
a_{ij}\right) \right] _{i,j=1}^n 
\]
is positive for each positive integer $n.$
\end{Definition}

In 1955, Stinespring [\textbf{18}] showed that a completely positive linear
map $\rho $ from a $C^{\ast }$-algebra $A$ to $L(H)$, the $C^{\ast }$%
-algebra of all bounded linear operators on a Hilbert space $H$, induces a
representation $\Phi _{\rho _{{}}}$ of $A$ on another Hilbert space $H_{\rho
_{{}}}$. Moreover, 
\[
\rho \left( a\right) =V_{\rho _{{}}}^{\ast }\Phi _{\rho _{{}}}(a)V_{\rho
_{{}}} 
\]
for all $a\in A$ and for some bounded linear operator $V_{\rho _{{}}}$ from $%
H$ to $H_{\rho _{{}}}$. \smallskip

In 1969, Arveson [\textbf{1}] proved a Radon-Nikodym type theorem which
gives a description of the order relation in the set of all completely
positive linear maps from $A$ to $L(H)$ in terms of the Stinespring
representation associated with each completely positive linear map. The gist
of the proof of this result is the fact that any bounded linear operator on
a Hilbert space is adjointable.

\smallskip Hilbert $C^{*}$-modules are generalizations of Hilbert spaces by
allowing the inner-product to take values in a $C^{*}$-algebra rather than
in the field of complex numbers.

\begin{Definition}
A pre-Hilbert $A$-module is a complex vector space$\ E$\ which is also a
right $A$-module, compatible with the complex algebra structure, equipped
with an $A$-valued inner product $\left\langle \cdot ,\cdot \right\rangle
:E\times E\rightarrow A\;$which is $\Bbb{C}$ -and $A$-linear in its second
variable and satisfies the following relations:

\begin{enumerate}
\item  $\ \ \left\langle \xi ,\eta \right\rangle ^{*}=\left\langle \eta ,\xi
\right\rangle \;\;$for every $\xi ,\eta \in E;$

\item  $\;\left\langle \xi ,\xi \right\rangle \geq 0\;\;$for every $\xi \in
E;$

\item  $\;\left\langle \xi ,\xi \right\rangle =0\;$\ if and only if $\xi =0.$

We say that $E\;$is a Hilbert $A$-module if $E\;$is complete with respect to
the topology determined by the norm $\left\| \cdot \right\| \;$given by $%
\left\| \xi \right\| =\sqrt{\left\| \left\langle \xi ,\xi \right\rangle
\right\| }.$
\end{enumerate}
\end{Definition}

\smallskip A \ $C^{*}$-algebra $A$ is a Hilbert $C^{*}$-module over $A$ with
the inner-product defined by $\left\langle a,b\right\rangle =a^{*}b$ for $a$
and $b$ in $A$.

Given two Hilbert $A$ -modules $E$ and $F$, the Banach space of all bounded
module homomorphisms from $E$ to $F$ is denoted by $\mathcal{B}_A(E,F)$. The
subset of $\mathcal{B}_A(E,F)$ consisting of all adjointable module
homomorphisms from $E$ to $F$ ( that is, $T\in \mathcal{B}_A(E,F)$ such that
there is $T^{*}$ $\in \mathcal{B}_A(F,E)$ satisfying $\left\langle \eta
,T\xi \right\rangle =\left\langle T^{*}\eta ,\xi \right\rangle $ for all $%
\xi \in E $ and for all $\eta \in F$ ) is denoted by $\mathcal{L}_A(E,F)$.
We will write $\mathcal{B}_A(E)$ for $\mathcal{B}_A(E,E)$ and $\mathcal{L}%
_A(E)$ for $\mathcal{L}_A(E,E) $.

In general, $\mathcal{L}_A(E,F)\neq \mathcal{B}_A(E,F)$. So the theory of
Hilbert $C^{*}$-modules is different from the theory of Hilbert spaces.

The Banach space $E^{\#}$ of all bounded module homomorphisms from $E$ to $A$
becomes a right $A$ -module with the action of $A$ on $E^{\#}$ defined by $%
\left( aT\right) \left( \xi \right) =a^{*}\left( T\xi \right) $ for $a\in A,$
$T\in E^{\#}$ and $\xi \in E$. We say that $E$ is self-dual if $E^{\#}=E$ as
right $A$ -modules.

If $E$ and $F$ are self-dual, then $\mathcal{B}_A(E,F)=\mathcal{L}_A(E,F)$ [%
\textbf{16, }Proposition 3.4].

\smallskip Suppose that \ $A$ is a $W^{*}$-algebra. Then the $A$ -valued
inner- product on $E$ extends to an $A$ -valued inner-product on $E^{\#}$
and in this way $E^{\#}$ becomes a self-dual Hilbert $A$ -module [\textbf{%
16, }Theorem 3.2]. Moreover, any bounded module homomorphism $T$ from\ $E$
to $F$ extends uniquely to a bounded homomorphism \ $\widetilde{T}$ from \ $%
E^{\#}$ to $F^{\#}$ [\textbf{16, }Proposition 3.6].

\smallskip \smallskip A representation of a $C^{*}$-algebra $A$ on a Hilbert 
$C^{*}$-module $E$ over a $C^{*}$-algebra $B$ is a $*$- morphism $\Phi $
from $A$ to $\mathcal{L}_B(E).$

Paschke [\textbf{16}] showed that a completely positive map from a unital $%
C^{*} $-algebra $A$ to another unital $C^{*}$-algebra $B$ induces a
representation of $A$ on a Hilbert $B$ -module which generalizes the GNS
construction and he extended the Arveson's results for completely positive
maps from a unital $C^{*}$-algebra $A$ to a $W^{*}$-algebra $B$.

In 1996, Tsui [\-\textbf{19}] proved a Radon-Nikodym type theorem for
completely positive maps between unital $C^{*}$-algebras and using this
theorem he obtained characterizations of the pure elements and the extreme
points in the set of all identity preserving completely positive maps from a
unital $C^{*}$-algebra $A$ to another unital $C^{*}$-algebra $B$ in terms of
a self-dual Hilbert module structure induced by each completely positive
map. To prove these facts he used the following construction.

\begin{Construction}
( [\textbf{15}, \textbf{16}, \textbf{19}]) Let $E$ be a Hilbert $C^{*}$%
-module over a $C^{*}$-algebra $B$. The algebraic tensor product $E\otimes _{%
\text{alg}}B^{**},$ where $B^{**\;}$ is the enveloping $W^{*}$ -algebra of $%
B,$ becomes a right $B^{**}$ -module if we define $\left( \xi \otimes
b\right) c=\xi \otimes bc,$ for $\xi \in E,$ and $b,c\in B^{**}.$

The map $\left[ \cdot ,\cdot \right] :$ $\left( E\otimes _{\text{alg}%
}B^{**}\right) \times \left( E\otimes _{\text{alg}}B^{**}\right) \rightarrow 
$ $B^{**}$ defined by 
\[
\left[ \sum\limits_{i=1}^n\xi _i\otimes b_i,\sum\limits_{j=1}^m\eta
_j\otimes c_j\right] =\sum\limits_{i=1}^n\sum\limits_{j=1}^mb_i^{*}\left%
\langle \xi _i,\eta _j\right\rangle c_j 
\]
is a $B^{**}$-valued inner-product on $E\otimes _{\text{alg}}B^{**}$ and the
quotient module $E\otimes _{\text{alg}}B^{**}/N_E$, where $N_E=\{\zeta \in
E\otimes _{\text{alg}}B^{**};\left[ \zeta ,\zeta \right] =0\}$, becomes a
pre-Hilbert $B^{**}$-module. The Hilbert $C^{*}$-module $\overline{E\otimes
_{\text{alg}}B^{**}/N_E}$ obtained by the completion of $E\otimes _{\text{alg%
}}B^{**}/N_E$ with respect to the norm induced by the inner product $\left[
\cdot ,\cdot \right] $ is called the extension of $E$ by the $C^{*}$
-algebra $B^{**}$. Moreover, $E$ can be regarded as a $B$-submodule of $%
\overline{E\otimes _{\text{alg}}B^{**}/N_E}$, since the map $\xi \mapsto \xi
\otimes 1+N_E$ from $E\;$to $\overline{E\otimes _{\text{alg}}B^{**}/N_E}$ is
an isometric inclusion.

The self-dual Hilbert $B^{**}$-module $\left( \overline{E\otimes _{\text{alg}%
}B^{**}/N_E}\right) ^{\#}$ is denoted by $\widetilde{E}$, and we can
consider $E$ as embedded in $\widetilde{E}$ without making distinction.
\end{Construction}

\smallskip Let $T\in \mathcal{B}_B(E,F).$ For $b_1,...,b_m\in B^{**}$ and $%
\xi _1,...,\xi _m$ in $E$ we denote by $b$ the element in $\left(
B^{**}\right) ^m$ whose components are $b_1,...,b_m,$ by $X$ the matrix in $%
M_n(B^{**})$ whose the $\left( i,j\right) $ -entry is $\left\langle \xi
_i,\xi _j\right\rangle $ and by $X_T$ the matrix in $M_n(B^{**})$ whose the $%
\left( i,j\right) $ -entry is $\left\langle T\xi _i,T\xi _j\right\rangle .$
By Lemma 4.2 in [\textbf{14}], $0\leq X_T\leq \left\| T\right\| X.$
Identifying $M_n(B^{**})$ with $L_{B^{**}}((B^{**})^n),$ we have 
\begin{eqnarray*}
\left[ \sum\limits_{i=1}^mT\xi _i\otimes b_i,\sum\limits_{i=1}^mT\xi
_i\otimes b_i\right] &=&\sum\limits_{i,j=1}^mb_i^{*}\left\langle T\xi
_i,T\xi _j\right\rangle b_j=\left\langle b,X_Tb\right\rangle \\
&\leq &\left\| T\right\| \left\langle b,Xb\right\rangle =\left\| T\right\| 
\left[ \sum\limits_{i=1}^m\xi _i\otimes b_i,\sum\limits_{i=1}^m\xi _i\otimes
b_i\right] .
\end{eqnarray*}
Therefore $T$ extends uniquely to a bounded module homomorphism $\widehat{T}$
from $\overline{E\otimes _{\text{alg}}B^{**}/N_E}$ to $\overline{F\otimes _{%
\text{alg}}B^{**}/N_F}$ such that 
\[
\widehat{T}\left( \sum\limits_{i=1}^m\xi _i\otimes b_i\right)
=\sum\limits_{i=1}^mT\xi _i\otimes b_i 
\]
and by Proposition 3.6 in [\textbf{16}], this extends uniquely to a bounded
module homomorphism $\widetilde{T}$ from $\widetilde{E}$ to $\widetilde{F}$
such that $\left\| T\right\| =\left\| \widetilde{T}\right\| .$

\begin{Remark}
\smallskip Any element $T\in \mathcal{B}_B(E,F)$ extends uniquely to an
element $\widetilde{T}\in \mathcal{B}_{B^{**}}(\widetilde{E},\widetilde{F})$
such that $\left\| T\right\| =\left\| \widetilde{T}\right\| $. Moreover, $%
\widetilde{TS}=\widetilde{T}\widetilde{S}$ for all $T\in \mathcal{B}_B(E,F)$
and $S\in \mathcal{B}_B(F,E),$ and if $T\in $ $\mathcal{L}(E,F),$ then $%
\widetilde{T^{*}}=\widetilde{T}^{*}$.
\end{Remark}

\smallskip

\begin{Remark}
Let $T$ $\in $ $\mathcal{B}_B(E,E^{\#})$. We extend\ \ \ $T$ to an element $%
\;\;\overline{T}$ $\in $ $\mathcal{B}_B$ $(\overline{\left( E\otimes _{\text{%
alg}}B^{**}\right) /N_E},\widetilde{E})$ by 
\[
\left[ \overline{T}\left( \sum\limits_{i=1}^n\xi _i\otimes b_i\right)
,\sum\limits_{j=1}^m\eta _j\otimes c_j\right] =\sum_{i=1}^n%
\sum_{j=1}^mb_i^{*}\left[ T\xi _i,\eta _j\right] c_j 
\]
and then extend it again to an element $\widetilde{T}\in \mathcal{B}_B(%
\widetilde{E})$ such that $\left\| T\right\| =\left\| \widetilde{T}\right\| $
[\-\textbf{16, }Proposition 3.6 ] .
\end{Remark}

\begin{Remark}
A representation $\Phi $ of a $C^{*}$-algebra $A$ on a Hilbert $C^{*}$%
-module $E$ over a $C^{*}$-algebra $B$ induces a representation $\widetilde{%
\Phi }$ of $A$ on $\widetilde{E}$ defined by $\widetilde{\Phi }\left(
a\right) =\widetilde{\Phi \left( a\right) }$ for all $a\in A.$
\end{Remark}

\begin{Remark}
Any completely positive linear map $\rho $ from $A$ to $B$ induces a
representation $\widetilde{\Phi _{\rho _{}}}$ of $A$ on a self -dual Hilbert 
$B^{**}$-module $\widetilde{E_{\rho _{}}}$ .
\end{Remark}

\smallskip Locally $C^{*}$-algebras are generalizations of $C^{*}$-algebras.
Instead of being given by a single norm, the topology on a locally $C^{*}$%
-algebra is defined by a directed family of $C^{*}$-seminorms. In fact a
locally $C^{*}$-algebra is an inverse limit of $C^{*}$-algebras.

\begin{Definition}
A locally $C^{*}$-algebra is a complete complex Hausdorff topological $*$
-algebra $A$ whose topology is determined by its continuous $C^{*}$
-seminorms in the sense that the net $\{a_i\}_{i\in I}$ converges to $0$ if
and only if the net $\{p(a_i)\}_{i\in I}$ converges to $0$ for all
continuous $C^{*}$-seminorm $p$ on $A.$
\end{Definition}

If $A$ is a locally $C^{*}$-algebra and $S(A)$ is the set of all continuous $%
C^{*}$-seminorms on $A$, then for each $p\in S(A),$ $A_p=A/\ker p$ is a $%
C^{*}$-algebra in the norm induced by $p$. The canonical map from $A$ onto $%
A_p$ is denoted by $\pi _p$ for each $p\in S(A)$. For $p,q\in S(A)$ with $%
q\leq p$ there is a unique morphism of $C^{*}$ -algebras $\pi _{pq}$ from $%
A_p$ onto $A_q$ such that $\pi _{pq}\left( \pi _p\left( a\right) \right)
=\pi _q\left( a\right) $ for all $a\in A.$ Moreover, $\{A_p;\pi
_{pq}\}_{p,q\in S(A),p\geq q}$ is an inverse system of $C^{*}$-algebras, and 
$A$ can be identified with $\lim\limits_{\stackunder{p}{\leftarrow }}A_p.$
Clearly, any $C^{*}$ -algebra is a locally $C^{*}$-algebra.

The terminology ''locally $C^{*}$-algebra'' is due to Inoue [\textbf{6}]. In
the literature, locally $C^{*}$-algebras have been given different name such
as $b^{*}$-algebras, $m$ -convex- $C^{*}$-algebras, $LMC^{*}$-algebras [%
\textbf{3}] or pro- $C^{*}$-algebras [\textbf{17}]. Such important concepts
as Hilbert $C^{*}$-modules, adjointable operators, (completely) positive
linear maps, (completely) multi-positive linear maps can be defined with
obvious modifications in the framework of locally $C^{*}$-algebras and many
results from the theory of $C^{*}$-algebras are still valid. The proofs are
not always straightforward. Thus, in [\textbf{3}] it is proved that a
continuous positive functional $\rho $ on a locally $C^{*}$-algebra $A$
induces a representation of $A$ on a Hilbert space $H$ which extends the GNS
construction, and moreover, the representation of $A$ induced by $\rho $ is
irreducible if and only if $\rho $ is pure. In [\textbf{2}], Bhatt and Karia
extend the Stinespring construction for completely positive linear maps from
a locally $C^{*}$-algebra $A$ to $L(H).$

If $A$ is a locally $C^{*}$ -algebra, then the set $M_n(A)$ of all $n\times
n $ matrices over $A$ with the algebraic operations and the topology
obtained by replying it as a direct sum of $n^2$ copies of $A$ is a locally $%
C^{*}$ -algebra.

\begin{Definition}
([\textbf{4}], [\textbf{9}]). A completely $n$ -positive map from a locally $%
C^{*}$-algebra $A$ to another locally $C^{*}$-algebra $B$ is an $n\times n$
matrix $\left[ \rho _{ij}\right] _{i,j=1}^n$whose elements are continuous
linear maps from $A$ to $B$ such that the map $\rho $ from $M_n(A)$ to $%
M_n(B)$ defined by 
\[
\rho \left( \left[ a_{ij}\right] _{i,j=1}^n\right) =\left[ \rho _{ij}\left(
a_{ij}\right) \right] _{i,j=1}^n 
\]
is completely positive.
\end{Definition}

\begin{Definition}
\smallskip ([\textbf{9}]). A completely $n$ -positive map $\left[ \rho _{ij}%
\right] _{i,j=1}^n$ from $A$ to $\mathcal{L}_B(E),$ where $E$ is a Hilbert
module over a $C^{*}$-algebra $B$ is strict if for some approximate unit $%
\{e_{\lambda _{}}\}_{\lambda \in \Lambda }$ for $A,$ the nets $\{\rho
_{ii}\left( e_{\lambda _{}}\right) \}_{\lambda \in \Lambda },$ $i\in
\{1,...,n\}$ are strictly Cauchy in $\mathcal{L}_B(E)$ ( that is, the nets $%
\{\rho _{ii}\left( e_{\lambda _{}}\right) \xi \}_{\lambda \in \Lambda },$ $%
i\in \{1,...,n\}$ are Cauchy in $E$ for each $\xi \in E)$.
\end{Definition}

\begin{Remark}
If $A$ is unital or $E$ is a Hilbert space, then any completely $n$-positive
map from $A$ to $\mathcal{L}(E)$ is strict.
\end{Remark}

\smallskip In [\textbf{9}], we extend the KSGNS (Kasparov, Stinespring,
Gel'fand, Naimark, Segal) construction for strict completely multi-positive
linear maps between locally $C^{*}$-algebras.

\begin{theorem}
(\textbf{\ } [\textbf{9}]). Let $A$ be a locally $C^{\ast }$-algebras, let $E
$ be a Hilbert module over a $C^{\ast }$-algebra $B$ and let $\rho =\left[
\rho _{ij}\right] _{i,j=1}^{n}$ be a strict completely $n$ -positive map
from $A$ to $\mathcal{L}_{B}(E).$

\begin{enumerate}
\item  There is a representation $\Phi _{\rho _{{}}}$ of $A$ on a Hilbert $B$
-module $E_{\rho _{{}}}$ and there are $n$ elements $V_{\rho ,1},...,V_{\rho
,n}$ in $\mathcal{L}_{B}\left( E,E_{\rho _{{}}}\right) $ such that

\begin{enumerate}
\item  $\rho _{ij}\left( a\right) =V_{\rho ,i}^{\ast }\Phi _{\rho
_{{}}}\left( a\right) V_{\rho ,j}$ for all $a\in A$ and for all $i,j\in
\{1,...,n\};$

\item  $\{\Phi _{\rho _{{}}}(a)V_{\rho ,i}x;a\in A,x\in E,1\leq i\leq n\}$
spans a dense subspace of $E_{\rho _{{}}}.$
\end{enumerate}

\item  If $\Phi $ is another representation of $A$ on a Hilbert $B$ -module $%
F$ and $W_{1},...,W_{n}$ are $n$ elements in $\mathcal{L}_{B}(E,F)$ such that

\begin{enumerate}
\item  $\rho _{ij}\left( a\right) =W_{i}^{\ast }\Phi \left( a\right) W_{j}$
for all $a\in A$ and for all $i,j\in \{1,...,n\};$

\item  $\{\Phi (a)W_{i}x;a\in A,x\in E,1\leq i\leq n\}$ spans a dense
subspace of $F,$

there is a unitary operator $U\in ${}$\mathcal{L}_{B}(E_{\rho _{{}}},F)$
such that

\begin{enumerate}
\item  $\Phi \left( a\right) U=U\Phi _{\rho _{{}}}(a)$ for all $a\in A;$ and

\item  $W_{i}=UV_{\rho ,i}$ for all $i\in \{1,...,n\}.$
\end{enumerate}
\end{enumerate}
\end{enumerate}
\end{theorem}

\smallskip The $n+2$ tuple $\left( \Phi _{\rho _{}},E_{\rho _{}},V_{\rho
,1},...,V_{\rho ,n}\right) $ is called the KSGNS construction associated
with $\rho .$

In [\textbf{10}], we prove a Radon-Nikodym type theorem for completely
multi- positive linear maps from a locally $C^{*}$-algebra $A$ to $L(H)$ and
we characterize the pure elements and the extreme points in the set of all
identity preserving completely multi-positive linear maps from $A$ to $L(H)$
in terms of the representation of $A$ induced by each completely
multi-positive linear map. Also, we determine a certain class of extreme
points in the set of all identity preserving completely positive linear maps
from $A$ to $M_n(L(H)).$ In this talk, we will extend the results from [%
\textbf{10}] for completely multi-positive linear maps from a locally $C^{*}$%
-algebra $A$ to a $C^{*}$-algebra $B$.

\section{The Radon-Nikodym theorem for completely $n$-positive linear maps}

Throughout this section, we assume that $A$ is a locally $C^{*}$-algebra, $B$
is a $C^{*}$ -algebra and $E$ is a Hilbert $C^{*}$-module over $B.$ We will
denote by $SCP_\infty ^n\left( A,\mathcal{L}_B(E)\right) $ the set of all
strict completely $n$ -positive linear maps from $A$ to $\mathcal{L}_B(E)$
and by $CP_\infty \left( A,\mathcal{L}_B(E)\right) $ the set of all
completely positive linear maps from $A$ to $\mathcal{L}_B(E).$

\begin{Proposition}
( [\textbf{4}],[\textbf{10}]) There is a bijection $\mathcal{S}$ from the
set $CP_\infty ^n(A,B)$ of all completely $n$-positive maps from $A$ to $B$
onto the set $CP_\infty (A,M_n(B))$ of all completely positive maps from $A$
to $M_n(B)$ defined by 
\[
\mathcal{S}\left( \left[ \rho _{ij}\right] _{i,j=1}^n\right) \left( a\right)
=\left[ \rho _{ij}(a)\right] _{i,j=1}^n\text{ for all }a\in A 
\]
which preserves the order relation.
\end{Proposition}

For an element $T\in \mathcal{L}_{B^{**}}(\widetilde{E})$ we denotes by $%
T|_E $ the restriction of the map $T$ on $E.$

Let $\rho \in SCP_\infty ^n\left( A,\mathcal{L}_B(E)\right) .$ We denote by $%
C\left( \rho \right) $ the $C^{*}$-subalgebra of $\mathcal{L}_{B^{**}}(%
\widetilde{E_{\rho _{}}})$ generated by $\{T\in \mathcal{L}_{B^{**}}(%
\widetilde{E_{\rho _{}}});T\widetilde{\Phi _{\rho _{}}}\left( a\right) =%
\widetilde{\Phi _{\rho _{}}}\left( a\right) T,$ $\left. \widetilde{V_{\rho
,j}}^{*}T\widetilde{\Phi }_{\rho _{}}\left( a\right) \widetilde{V_{\rho ,i}}%
\right| _E $ $\in \mathcal{L}_B(E)$ for all $a\in A$ and for all $i,j\in
\{1,...,n\}\}.$

\begin{Remark}
If $T$ is an element in $C\left( \rho \right) $, then $T|_{E_{\rho _{}}}\in 
\mathcal{B}_B(E_{\rho _{}},E_{\rho _{}}^{\#})$, since 
\[
\left\langle T\Phi _{\rho _{}}\left( a\right) V_{\rho ,j}\xi ,\Phi _{\rho
_{}}\left( b\right) V_{\rho ,i}\eta \right\rangle =\left\langle T\widetilde{%
\Phi }_{\rho _{}}\left( a\right) \widetilde{V_{\rho ,j}}\xi ,\widetilde{\Phi 
}_{\rho _{}}\left( b\right) \widetilde{V_{\rho ,i}}\eta \right\rangle \in B 
\]
for all $a,b\in A,$ for all $\xi ,\eta \in E$ and for all $i,j\in
\{1,...,n\} $ and since $\{\Phi _{\rho _{}}\left( a\right) V_{\rho ,i}\xi ;$ 
$a\in A,\xi \in E,1\leq i\leq n\}$ spans a dense submodule of $E_{\rho _{}}.$
\end{Remark}

\begin{Lemma}
Let $T\in C\left( \rho \right) .$ If $T$ is positive, then the map $\rho _T$
from $M_n(A)$ to $M_n(\mathcal{L}_B(E))$ defined by 
\[
\rho _T\left( \left[ a_{ij}\right] _{i,j=1}^n\right) =\left[ \left. 
\widetilde{V_{\rho ,i}}^{*}T\widetilde{\Phi _{\rho _{}}}\left( a_{ij}\right) 
\widetilde{V_{\rho ,j}}\right| _E\right] _{i,j=1}^n 
\]
is a strict completely $n$ -positive linear map from $A$ to $\mathcal{L}%
_B(E).$
\end{Lemma}

\proof%
%
It is not difficult to see that $\rho _T$ is an $n\times n$ matrix of
continuous linear maps from $A$ to $\mathcal{L}_B(E),$ the $\left(
i,j\right) $-entry of the matrix $\rho _T$ is the linear map $\left( \rho
_T\right) _{ij} $ from $A$ to $\mathcal{L}_B(E)$ defined by 
\[
\left( \rho _T\right) _{ij}\left( a\right) =\left. \widetilde{V_{\rho ,i}}%
^{*}T\widetilde{\Phi _{\rho _{}}}\left( a\right) \widetilde{V_{\rho ,j}}%
\right| _E. 
\]
Also it is not difficult to check that for all $a_1,...,a_m\in A$ and for
all $T_1,...,T_m\in M_n(\mathcal{L}_B(E)),$ we have 
\begin{eqnarray*}
\sum\limits_{k,l=1}^mT_l^{*}\mathcal{S}\left( \rho _T\right) \left(
a_l^{*}a_k\right) T_k &=&\left. \left( \sum\limits_{k,l=1}^m\widetilde{T_l}%
^{*}\widetilde{\mathcal{S}\left( \rho _T\right) }\left( a_l^{*}a_k\right) 
\widetilde{T_k}\right) \right| _E \\
&=&\left. \left( \left( \sum\limits_{l=1}^mM_{T^{\frac 12}}(a_l)V\widetilde{%
T_l}\right) ^{*}\left( \sum\limits_{l=1}^mM_{T^{\frac 12}}(a_l)V\widetilde{%
T_l}\right) \right) \right| _E,
\end{eqnarray*}
where $M_{T^{\frac 12}}(a)=\left[ 
\begin{array}{ccc}
T^{\frac 12}\widetilde{\Phi _{\rho _{}}}\left( a\right) & \cdot \cdot \cdot
& T^{\frac 12}\widetilde{\Phi _{\rho _{}}}\left( a\right) \\ 
0 & \cdot \cdot \cdot & 0 \\ 
\cdot & \cdot \cdot \cdot & \cdot \\ 
0 & \cdot \cdot \cdot & 0
\end{array}
\right] $ and $V=\left[ 
\begin{array}{ccc}
\widetilde{V_{\rho ,1}} & \cdot \cdot \cdot & 0 \\ 
\cdot & \cdot \cdot \cdot & \cdot \\ 
0 & \cdot \cdot \cdot & \widetilde{V_{\rho ,n}}
\end{array}
\right] .$ From this fact we conclude that $\mathcal{S}\left( \rho _T\right)
\in CP_\infty (A,M_n(\mathcal{L}_B(E)))$ and by Proposition 2.1, $\rho _T\in
CP_\infty ^n(A,\mathcal{L}_B(E)).$ To show that $\rho _T\in SCP_\infty ^n(A,%
\mathcal{L}_B(E)),$ let $\{e_{\lambda _{}}\}_{\lambda \in \Lambda }$ be an
approximate unit for $A$, $\xi \in E$ and $i\in \{1,...,n\}.$ Then 
\begin{eqnarray*}
\left\| \left( \rho _T\right) _{ii}(e_{\lambda _{}})\xi -\left( \rho
_T\right) _{ii}(e_{\mu _{}})\xi \right\| &=&\left\| \widetilde{V_{\rho ,i}}%
^{*}T\left( \widetilde{\Phi _{\rho _{}}}(e_{\lambda _{}})-\widetilde{\Phi
_{\rho _{}}}(e_{\mu _{}})\right) \widetilde{V_{\rho ,i}}\xi \right\| \\
&\leq &\left\| \widetilde{V_{\rho ,i}}^{*}T\right\| \left\| \left( \Phi
_{\rho _{}}(e_{\lambda _{}})-\Phi _{\rho _{}}(e_{\mu _{}})\right) V_{\rho
,i}\xi \right\| ,
\end{eqnarray*}
and since $\{\Phi _{\rho _{}}(e_{\lambda _{}})V_{\rho ,i}\xi \}_{\lambda \in
\Lambda }$ is a Cauchy net in $E,$ the net $\{\left( \rho _T\right)
_{ii}(e_{\lambda _{}})\}_{\lambda \in \Lambda }$ is strictly Cauchy.
Therefore $\rho _T\in SCP_\infty ^n(A,\mathcal{L}_B(E)).$%
\endproof%
%

\begin{Remark}
\begin{enumerate}
\item  If $I_{\widetilde{E_{\rho _{}}}}$ is the identity map on $\widetilde{%
E_{\rho _{}}},$ then $\rho _{I_{\widetilde{E_{\rho _{}}}}}=\rho .$

\item  If $T_1$ and $T_2$ are two positive elements in $C\left( \rho \right)
,$ then $\rho _{T_1+T_2}=\rho _{T_1}+\rho _{T_2}.$

\item  If $T$ is a positive element in $C\left( \rho \right) $ and $\alpha $
is a positive number, then $\rho _{\alpha T}=\alpha \rho _T.$
\end{enumerate}
\end{Remark}

\begin{Remark}
Let $T_1$ and $T_2$ be two positive elements in $C\left( \rho \right) .$ If $%
T_1\leq T_2,$ then, since 
\begin{eqnarray*}
\left( \rho _{T_2}-\rho _{T_1}\right) \left( \left[ a_{ij}\right]
_{i,j=1}^n\right) &=&\left[ \left. \widetilde{V_{\rho ,i}}^{*}\left(
T_2-T_1\right) \widetilde{\Phi }_{\rho _{}}\left( a_{ij}\right) \widetilde{%
V_{\rho ,j}}\right| _E\right] _{i,j=1}^n \\
&=&\rho _{T_2-T_1}\left( \left[ a_{ij}\right] _{i,j=1}^n\right)
\end{eqnarray*}
for all $\left[ a_{ij}\right] _{i,j=1}^n\in M_n(A),$ $\rho _{T_1}\leq \rho
_{T_2}.$
\end{Remark}

Let $\rho \in SCP_\infty ^n(A,\mathcal{L}_B(E)).$ We denote by $\left[
0,\rho \right] $ the set of all strict completely $n$ -positive linear maps $%
\theta $ from $A$ to $\mathcal{L}_B(E)$ such that $\theta \leq \rho $ ( that
is, $\rho -\theta \in SCP_\infty ^n(A,\mathcal{L}_B(E))$ ) and by $\left[ 0,I%
\right] _{\rho _{}}$ the set of all elements $T$ in $C\left( \rho \right) $
such that $0\leq T\leq I_{\widetilde{E_{\rho _{}}}}.$

\begin{Theorem}
The map $T$ $\rightarrow \rho _T$ from $\left[ 0,I\right] _{\rho _{}}$ to $%
\left[ 0,\rho \right] $ is an affine order isomorphism.
\end{Theorem}

\proof%
%
By Lemma 2.3 and Remarks 2.4 and 2.5, the map $T$ $\rightarrow \rho _T$ from 
$\left[ 0,I\right] _{\rho _{}}$ to $\left[ 0,\rho \right] $ is well-defined
and moreover, it is affine. To show that this map is injective, let $T\in %
\left[ 0,I\right] _{\rho _{}}$ such that $\rho _T=0.$ Then $\left. 
\widetilde{V_{\rho ,i}}^{*}T\widetilde{\Phi }_{\rho _{}}\left( a\right) 
\widetilde{V_{\rho ,j}}\right| _E=0$ for all $a\in A$ and for all $i,j\in
\{1,2,...,n\},$ and so 
\[
\left\langle T\Phi _{\rho _{}}(a)V_{\rho ,j}\xi ,\Phi _{\rho _{}}\left(
b\right) V_{\rho ,j}\eta \right\rangle =0 
\]
for all $a,b\in A$ for all $\xi ,\eta \in E$ and for all $i,j\in
\{1,...,n\}. $ Taking into account that $\{\Phi _{\rho _{}}(a)V_{\rho ,i}\xi
;a\in A,\xi \in E,1\leq i\leq n\}$ spans a dense submodule of $E_{\rho _{}},$
from these facts, Remarks 2.2 and 1.5 and we conclude that $T=0.$

Let $\theta \in \left[ 0,\rho \right] .$ In the same way as in the proof of
Lemma 3.4 in [\textbf{10}], we show that there is a bounded linear map $W$
from $E_{\rho _{}}$ to $E_{\theta _{}}$ such that 
\[
W\left( \Phi _{\rho _{}}\left( a\right) V_{\rho ,i}\xi \right) =\Phi
_{\theta _{}}\left( a\right) V_{\theta ,i}\xi . 
\]
It is not difficult to check that $W$ is a bounded module homomorphism such
that $\left\| W\right\| \leq 1,$ $WV_{\rho ,i}=V_{\theta ,i}$ for all $i\in
\{1,...,n\}$ and $W\Phi _{\rho _{}}\left( a\right) =\Phi _{\theta _{}}\left(
a\right) W$ for all $a\in A.$ If $\widetilde{W}$ is the unique extension of $%
W$ to a bounded module morphism from $\widetilde{E_{\rho _{}}}$ to $%
\widetilde{E_{\theta _{}}}$ with $\left\| W\right\| =\left\| \widetilde{W}%
\right\| ,$ then clearly $0\leq \widetilde{W}^{*}\widetilde{W}\leq I_{%
\widetilde{E_{\rho _{}}}}.$ Moreover, it is easy to check that $\widetilde{W}%
^{*}\widetilde{W}\widetilde{\Phi _{\rho _{}}}\left( a\right) =\widetilde{%
\Phi _{\rho _{}}}\left( a\right) \widetilde{W}^{*}\widetilde{W}$ for all $%
a\in A,$ and since 
\[
\widetilde{V_{\rho ,i}}^{*}\widetilde{W}^{*}\widetilde{W}\widetilde{\Phi
_{\rho _{}}}\left( a\right) \widetilde{V_{\rho ,j}}=\widetilde{V_{\rho ,i}}%
^{*}\widetilde{W}^{*}\widetilde{\Phi _{\theta _{}}}\left( a\right) 
\widetilde{W}\widetilde{V_{\rho ,j}}=\widetilde{V_{\theta ,i}}^{*}\widetilde{%
\Phi _{\theta _{}}}\left( a\right) \widetilde{V_{\theta ,j}} 
\]
for all $a\in A$ and for all $i,j\in \{1,...,n\},$ $\widetilde{W}^{*}%
\widetilde{W}$ $\in \left[ 0,I\right] _{\rho _{}}.$ Let $T=\widetilde{W}^{*}%
\widetilde{W}.$ Then clearly, $\theta =\rho _T$ and thus the map $T$ $%
\rightarrow \rho _T$ from $\left[ 0,I\right] _{\rho _{}}$ to $\left[ 0,\rho %
\right] $ is surjective. Therefore the map $T$ $\rightarrow \rho _T$ is an
affine isomorphism from $\left[ 0,I\right] _{\rho _{}}$ onto $\left[ 0,\rho %
\right] $ which preserve the order relation.%
\endproof%
%

\section{Applications of the Radon-Nikodym theorem}

Let $A$ be a locally $C^{*}$ -algebra, let $B$ be a $C^{*}$ -algebra and let 
$E$ be a Hilbert $C^{*}$-module over $B.$ A strict completely $n$-positive
linear map $\rho $ from $A$ to $\mathcal{L}_B(E)$ is said to be pure if for
every strict completely $n$-positive linear map $\theta $ from $A$ to $%
\mathcal{L}_B(E)$ with $\theta \leq \rho ,$ there is a positive number $%
\alpha $ such that $\theta =\alpha \rho .$

\begin{Proposition}
Let $\rho \in SCP_\infty ^n(A,\mathcal{L}_B(E)).$ Then $\rho $ is pure if
and only if $\left[ 0,I\right] _{\rho _{}}=\{\alpha I_{\widetilde{E_{\rho
_{}}}};0\leq \alpha \leq 1\}.$
\end{Proposition}

\proof%
%
First we suppose that $\rho $ is pure. Let $T\in \left[ 0,I\right] _{\rho
_{}}.$ By Theorem 2.6, $\rho _T\in \left[ 0,\rho \right] ,$ and since $\rho $
is pure, $\rho _T=\alpha \rho $ for some positive number. From this fact,
Remark 2.4 and Theorem 2.6 we deduce that $T=\alpha I_{\widetilde{E_{\rho
_{}}}}$ for some $0\leq \alpha \leq 1.$

Conversely, suppose that $\left[ 0,I\right] _{\rho _{}}=\{\alpha I_{%
\widetilde{E_{\rho _{}}}};0\leq \alpha \leq 1\}.$ Let $\theta \in SCP_\infty
^n(A,$ $\mathcal{L}_B(E))$ such that $\theta \leq \rho .$ By Theorem 2.6, $%
\theta =\rho _T$ for some $T\in $ $\left[ 0,I\right] _{\rho _{}},$ and since 
$T=\alpha I_{\widetilde{E_{\rho _{}}}}$ for some positive number $\alpha ,$ $%
\theta =\alpha \rho $ and the proposition is proved.%
\endproof%
%

\begin{Corollary}
A strict completely $n$ -positive linear map $\rho $ from $A$ to $\mathcal{L}%
_B(E)$ is pure if and only if $C\left( \rho \right) $ consisting of the
scalar multipliers of $I_{\widetilde{E_{\rho _{}}}}.$
\end{Corollary}

We say that two strict completely $n$-positive linear maps $\rho $ and $%
\theta $ from $A$ to $\mathcal{L}_B(E)\;$are unitarily equivalent if the
representations of $A$ induced by $\rho $ respectively $\theta $ are
unitarily equivalent.

The following proposition is a generalization of Proposition 4.3 in [\textbf{%
10}].

\begin{Proposition}
Let $A$ be a unital locally $C^{*}$-algebra, let $B$ be a $C^{*}$-algebra,
let $E$ be a Hilbert $B$ -module and let $\rho =\left[ \rho _{ij}\right]
_{i,j=1}^n\in CP_\infty ^n\left( A,\mathcal{L}_B(E)\right) .$ If $\rho
_{ii}, $ $i\in \{1,...,n\}$ are unitarily equivalent pure unital completely
positive linear maps from $A$ to $\mathcal{L}_B(E)$ and for all $i,j\in
\{1,...,n\}$ with $i\neq j$ there is a unitary element $u_{ij}$ in $A$ such
that $\rho _{ij}\left( u_{ij}\right) =I_E,$ then $\rho $ is pure.
\end{Proposition}

\proof%
%
Let $i,j\in \{1,...,n\}$ with $i\neq j.$ From 
\[
\left\| \Phi _{\rho _{}}\left( u_{ij}\right) V_{\rho ,j}-V_{\rho ,i}\right\|
^2=\left\| V_{\rho ,j}^{*}V_{\rho ,j}-\rho _{ij}(u_{ij})-\left( \rho
_{ij}(u_{ij})\right) ^{*}+V_{\rho ,i}^{*}V_{\rho ,i}\right\| =0 
\]
we deduce that $\Phi _{\rho _{}}\left( u_{ij}\right) V_{\rho ,j}=V_{\rho
,i}. $ Therefore the sets $\{\Phi _{\rho _{}}(a)V_{\rho ,i}\xi ;a\in A,\xi
\in E\} $ and $\{\Phi _{\rho _{}}(a)V_{\rho ,j}\xi ;a\in A,\xi \in E\}$
generate the same Hilbert submodule of $E_{\rho _{}},$ and since $E_{\rho
_{}}$ is generated by $\{\Phi _{\rho _{}}(a)V_{\rho ,i}\xi ;a\in A,\xi \in
E,1\leq i\leq n\},$ this coincides with $E_{\rho _{}}.$

Let $i\in \{1,...,n\}$ and let $(\Phi _{\rho _{ii}},E_{\rho _{ii}},V_{\rho
_{ii}})$ be the KSGNS\ construction associated with $\rho _{ii}.$ We will
show that the representations $\Phi _{\rho _{}}$ and $\Phi _{\rho _{ii}}$ of 
$A$ are unitarily equivalent. Since $\{\Phi _{\rho _{}}(a)V_{\rho ,i}\xi
;a\in A,\xi \in E\}$ spans a dense submodule of $E_{\rho _{}},$ $\{\Phi
_{\rho _{ii}}(a)V_{\rho _{ii}}\xi ;a\in A,\xi \in E\}$ spans a dense
submodule of $E_{\rho _{ii}}$ and 
\begin{eqnarray*}
\left\langle \Phi _{\rho _{}}(a)V_{\rho ,i}\xi ,\Phi _{\rho _{}}(b)V_{\rho
,i}\eta \right\rangle &=&\left\langle \rho _{ii}\left( b^{*}a\right) \xi
,\eta \right\rangle \\
&=&\left\langle V_{\rho _{ii}}^{*}\Phi _{\rho _{ii}}(b^{*}a)V_{\rho
_{ii}}\xi ,\eta \right\rangle \\
&=&\left\langle \Phi _{\rho _{ii}}(a)V_{\rho _{ii}}\xi ,\Phi _{\rho
_{ii}}(b)V_{\rho _{ii}}\eta \right\rangle
\end{eqnarray*}
for all $a,b\in A$ and for all $\xi ,\eta \in E,$ there is a unitary
operator $U_i$ from $E_{\rho _{ii}}$ to $E_{\rho _{}}$ such that 
\[
U_i\left( \Phi _{\rho _{ii}}(a)V_{\rho _{ii}}\xi \right) =\Phi _{\rho
_{}}(a)V_{\rho ,i}\xi 
\]
[\textbf{14, }Theorem 3.5]. Moreover, $U_i\Phi _{\rho _{ii}}(a)=$ $\Phi
_{\rho _{}}(a)$ $U_i$ for all $a\in A.$ Then $\widetilde{U_i},$ the unique
extension of $U_i$ to a bounded module homomorphism from $\widetilde{E_{\rho
_{ii}}}$ to $\widetilde{E_{\rho _{}}}$ with $\left\| U_i\right\| =\left\| 
\widetilde{U_i}\right\| ,$ is a unitary element in $\mathcal{L}_{B^{**}}(%
\widetilde{E_{\rho _{ii}}},\widetilde{E_{\rho _{}}}).$

Let $T\in \left[ 0,I\right] _{\rho _{}}.$ Then $\widetilde{U_i}^{*}T%
\widetilde{U_i}\in \left[ 0,I\right] _{\rho _{ii}},$ and since $\rho _{ii}$
is pure, by Proposition 3.1, $\widetilde{U_i}^{*}T\widetilde{U_i}=\alpha I_{%
\widetilde{E_{\rho _{ii}}}}$ for some positive number $\alpha .$
Consequently, $T=\alpha I_{\widetilde{E_{\rho _{}}}}$ and so $\rho $ is pure.%
\endproof%
%

In the following corollary we determine a class of extreme points in the set
of all identity preserving completely positive linear maps from a unital
locally $C^{*}$-algebra $A$ to the $C^{*}$-algebra $M_n(B)$ of all $n\times
n $ matrices with elements in the unital $C^{*}$-algebra $B.$ This is a
generalization of Corollaries 2.7 in [\textbf{11}] and 4.5 in [\textbf{10}].

\begin{Corollary}
Let $A$ be a unital locally $C^{*}$-algebra, let $B$ be a unital $C^{*}$%
-algebra and let $\rho =\left[ \rho _{ij}\right] _{i,j=1}^n\in CP_\infty
^n\left( A,B\right) .$ If $\rho _{ii},$ $i\in \{1,...,n\}$ are unitarily
equivalent pure unital completely positive linear maps from $A$ to $B$ and
for all $i,j\in \{1,...,n\}$ with $i\neq j,$ $\rho _{ij}\left( 1\right) =0$
and there is a unitary element $u_{ij}$ in $A$ such that $\rho _{ij}\left(
u_{ij}\right) =1,$ then the map $\varphi $ from $A$ to $M_n(B)$ defined by $%
\varphi \left( a\right) =\left[ \rho _{ij}\left( a\right) \right] _{i,j=1}^n$
is an extreme point in the set of all identity preserving completely
positive linear maps from $A$ to $M_n(B).$
\end{Corollary}

\proof%
%
Let $\varphi _1$ and $\varphi _2$ be two identity preserving completely
positive linear maps from $A$ to $M_n(B)$ and let $\alpha \in (0,1)$ such
that $\alpha \varphi _1+\left( 1-\alpha \right) \varphi _2=\varphi .$ Then $%
\alpha \mathcal{S}^{-1}(\varphi _1)+\left( 1-\alpha \right) \mathcal{S}%
^{-1}(\varphi _2)=\rho .$ From this relation and Propositions 3.3 and 3.1,
we conclude that $\alpha \mathcal{S}^{-1}(\varphi _1)=\beta _1\rho $ for
some positive number $\beta _1$ and $\left( 1-\alpha \right) \mathcal{S}%
^{-1}(\varphi _2)=\beta _2\rho $ for some positive number $\beta _2.$
Consequently, $\alpha \varphi _1=\beta _1\varphi $ and $\left( 1-\alpha
\right) \varphi _2=\beta _2\varphi .$ From these facts, since $\varphi
_1(1)=\varphi _2\left( 1\right) =\varphi \left( 1\right) =I_n,$ where $I_n$
is the unity matrix in $M_n(B),$ we deduce that $\alpha =\beta _1$ and $%
1-\alpha =\beta _2.$ Therefore $\varphi _1=\varphi _2=\varphi ,$ and so $%
\varphi $ is an extreme point in the set of all identity preserving
completely positive linear maps from $A$ to $M_n(B).$%
\endproof%
%

Let $A$ be a unital locally $C^{*}$-algebra, let $B$ be a $C^{*}$-algebra
and let $E$ be a Hilbert $B$ -module. We denote by $CP_\infty ^n(A,\mathcal{L%
}_B(E),I)$ the set of all completely $n$-positive linear maps $\rho =\left[
\rho _{ij}\right] _{i,j=1}^n$ from $A$ to $\mathcal{L}_B(E)$ such that $\rho
_{ii}(1)=I_E$ for all $i\in \{1,...,n\}$ and $\rho _{ij}\left( 1\right) =0$
for all $i,j\in \{1,..,n\}$ with $i\neq j.$

The following theorem is a generalization of Theorems 3.8 in [\textbf{18}]
and 4.6 in [\textbf{10}].

\begin{Theorem}
Let $\rho \in CP_\infty ^n(A,\mathcal{L}_B(E),I).$ Then $\rho $ is an
extreme point in the set $CP_\infty ^n(A,\mathcal{L}_B(E),I)$ if and only if
the map $T\rightarrow \left[ \widetilde{V_{\rho ,i}}^{*}T\widetilde{V_{\rho
,j}}\right] _{i,j=1}^n$ from $C\left( \rho \right) $ to $M_n(\mathcal{L}%
_{B^{**}}(\widetilde{E}))$ is injective.
\end{Theorem}

\proof%
%
Suppose that $\rho $ is an extreme point in the set $CP_\infty ^n(A,\mathcal{%
L}_B(E),I)$ and $T$ is an element in $C\left( \rho \right) $ such that $%
\widetilde{V_{\rho ,j}}^{*}T\widetilde{V_{\rho ,i}}=0$ for all $i,j\in
\{1,...,n\}.$ Since $\widetilde{V_{\rho ,j}}^{*}T^{*}\widetilde{V_{\rho ,i}}%
=\left( \widetilde{V_{\rho ,i}}^{*}T\widetilde{V_{\rho ,j}}\right) ^{*}$ for
all $i,j\in \{1,...,n\},$ we can suppose that $T=T^{*}.$ It is not difficult
to check that there are two positive numbers $\alpha $ and $\beta $ such
that $\frac 14I_{\widetilde{E_{\rho _{}}}}\leq \alpha T+\beta I_{\widetilde{%
E_{\rho _{}}}}\leq \frac 34I_{\widetilde{E_{\rho _{}}}}.$ Moreover, $\beta
\in \left( 0,1\right) .$ Let $T_1=\frac \alpha \beta T+I_{\widetilde{E_{\rho
_{}}}}$ and $T_2=I_{\widetilde{E_{\rho _{}}}}-\frac \alpha {1-\beta }T.$
Clearly, $T_k,$ $k\in \{1,2\}$ are two positive elements in $C\left( \rho
\right) .$ Then $\rho _{T_k}\in CP_\infty ^n(A,\mathcal{L}_B(E)),$ $k\in
\{1,2\} $ and since 
\[
\left( \rho _{T_k}\right) _{ij}(1)=\left. \widetilde{V_{\rho ,i}}^{*}T_k%
\widetilde{V_{\rho ,j}}\right| _E=\left. \widetilde{V_{\rho ,i}}^{*}%
\widetilde{V_{\rho ,j}}\right| _E=V_{\rho ,i}^{*}V_{\rho ,j}=\left\{ 
\begin{array}{ccc}
I_E & \text{if } & i=j \\ 
0 & \text{if } & i\neq j
\end{array}
\right. 
\]
for all $i,j\in \{1,...,n\}$ with $i\neq j,$ and $k\in \{1,2\},$ $\rho
_{T_k}\in CP_\infty ^n(A,\mathcal{L}_B(E),I)$ for each $k\in \{1,2\}.$ A
simple calculus shows that $\beta \rho _{T_1}+\left( 1-\beta \right) \rho
_{T_2}=\rho ,$ and since $\rho $ is an extreme point, $\rho _{T_1}=\rho
_{T_2}=\rho .$ But $\rho _{T_1}=\frac \alpha \beta \rho _T+\rho $ and $\rho
_{T_2}=\rho -\frac \alpha {1-\beta }\rho _T.$ Therefore $\rho _T=0$ and by
Theorem 2.6, $T=0.$

Conversely, suppose that the map\ $T\rightarrow \left[ \widetilde{V_{\rho ,i}%
}^{*}T\widetilde{V_{\rho ,j}}\right] _{i,j=1}^n$ \ from $C\left( \rho
\right) $ to $M_n(\mathcal{L}_{B^{**}}(\widetilde{E}))$ is injective. Let $%
\theta ,\sigma \in CP_\infty ^n(A,\mathcal{L}_B(E),I)$ and $\alpha \in
\left( 0,1\right) $ such that $\alpha \theta +\left( 1-\alpha \right) \sigma
=\rho . $ By Theorem 2.6, there are two elements $T_1$ and $T_2$ in $%
[0,I]_{\rho _{}}\subseteq C\left( \rho \right) ,$ such that $\alpha \theta
=\rho _{T_1}$ and $\left( 1-\alpha \right) \sigma =\rho _{T_2}.$ Then 
\[
\left. \widetilde{V_{\rho ,i}}^{*}T_1\widetilde{V_{\rho ,j}}\right|
_E=\left( \rho _{T_1}\right) _{ij}(1)=\alpha \theta _{ij}(1)=\left\{ 
\begin{array}{ccc}
\alpha I_E & \text{if } & i=j \\ 
0 & \text{if } & i\neq j
\end{array}
\right. 
\]
and 
\[
\left. \widetilde{V_{\rho ,i}}^{*}T_2\widetilde{V_{\rho ,j}}\right|
_E=\left( \rho _{T_2}\right) _{ij}(1)=\left( 1-\alpha \right) \sigma
_{ij}(1)=\left\{ 
\begin{array}{ccc}
\left( 1-\alpha \right) I_E & \text{if } & i=j \\ 
0 & \text{if } & i\neq j
\end{array}
\right. . 
\]
Therefore 
\[
\left. \widetilde{V_{\rho ,i}}^{*}\left( T_1-\alpha I_{\widetilde{E}_\rho
}\right) \widetilde{V_{\rho ,j}}\right| _E=\left. \widetilde{V_{\rho ,i}}%
^{*}\left( T_2-\left( 1-\alpha \right) I_{\widetilde{E}_\rho }\right) 
\widetilde{V_{\rho ,j}}\right| _E=0 
\]
from all $i,j\in \{1,...,n\}.$ From these facts, Remark 1.4 and taking into
account that $\left. \widetilde{V_{\rho ,i}}^{*}T_1\widetilde{V_{\rho ,j}}%
\right| _E$ and $\left. \widetilde{V_{\rho ,i}}^{*}T_2\widetilde{V_{\rho ,j}}%
\right| _E$ are elements in $\mathcal{L}_B(E)$ for all $i,j\in \{1,...,n\}$,
we conclude that 
\[
\widetilde{V_{\rho ,i}}^{*}\left( T_1-\alpha I_{\widetilde{E}_\rho }\right) 
\widetilde{V_{\rho ,j}}=\widetilde{V_{\rho ,i}}^{*}\left( T_2-\left(
1-\alpha \right) I_{\widetilde{E}_\rho }\right) \widetilde{V_{\rho ,j}}=0 
\]
for all $i,j\in \{1,...,n\}.$ Hence $T_1=\alpha I_{\widetilde{E_{\rho _{}}}}$
and $T_2=\left( 1-\alpha \right) I_{\widetilde{E_{\rho _{}}}}$ and
consequently $\theta =\sigma =\rho .$%
\endproof%
%

\textsc{Department of Mathematics}

\textsc{Faculty of Chemistry }

\textsc{University of Bucharest}

\textit{\ }\textsc{Bd. Regina Elisabeta nr.4-12}

\textsc{\ Bucharest, Romania\ \ }

\ \ \ \ \ \ \ \ mjoita@fmi.unibuc.ro

\end{document}